\documentclass[11pt]{amsart}
\usepackage{latexsym, amssymb,xypic,lscape,epsfig,setspace,hyperref}

\theoremstyle{plain}
\newtheorem{theo}{Theorem}[section]

\newtheorem{prop}[theo]{Proposition}

\theoremstyle{definition}
\newtheorem{example}[theo]{Example}
\newtheorem{definition}[theo]{Definition}

\newenvironment{renumerate}
{
\begin{enumerate}}
{\end{enumerate}
}

\newenvironment{remark}
{\vskip6pt \noindent {\it Remark:}} {\vskip6pt}

\newenvironment{remarks}%
{\vskip6pt%
\noindent%
{\it Remarks}: %
\begin{renumerate}}%
{\end{renumerate}\vskip6pt}
{\begin{example}\label{#1}}%
{\end{example}}

\newcommand{\map}{\Psi}

\renewcommand{\frak}[1]{\text{$\mathfrak{#1}$}}
\newcommand{\J}{\text{$\mathcal{J}$}}
\newcommand{\JJ}{\text{$\mathcal{J}$}}
\newcommand{\G}{\text{$\mathcal{G}$}}

\newcommand{\ga}{\text{$a$}}

\renewcommand{\gg}{\text{$u$}}
\newcommand{\gh}{\text{$w$}}

\newcommand{\Id}{\mathrm{Id}}

\newcommand{\delbar}{\text{$\overline{\partial}$}}

\newcommand{\mc}[1]{\text{$\mathcal{#1}$}}

\newcommand{\into}{\longrightarrow}
\newcommand{\noqed}{\let\qed\relax}

\newcommand{\Gg}{\mathfrak{g}}
\newcommand{\Aa}{\mathfrak{a}}
\newcommand{\Hh}{\mathfrak{h}}

\newcommand{\sym}{\mathbf{sym}}
\newcommand{\Ad}{\mathrm{ad}}

\newcommand{\gcs}{generalized complex structure}
\newcommand{\gacs}{generalized almost complex structure}
\newcommand{\gcss}{generalized complex structures}

\newcommand{\gk}{generalized K\"ahler}
\newcommand{\gks}{generalized K\"ahler structure}

\newcommand{\ghk}{generalized hyper-K\"ahler}
\newcommand{\ghks}{generalized hyper-K\"ahler structure}

\newcommand{\Ann}{\mathrm{Ann}}

\newcommand{\Kperp}{\text{$K^{\perp}$}}

\newcommand{\lra}{\longrightarrow}

\newcommand{\IP}[1]{\langle #1\rangle}

\newcommand{\E}{{E}}
\newcommand{\ad}{\mathrm{ad}}

\newcommand{\Cour}[1]{[\![#1]\!]}

\begin{document}

\title{Generalized K\"ahler and hyper-K\"ahler quotients}
\author{Henrique Bursztyn} \address{Instituto de Matem\'atica Pura
e Aplicada, Estrada Dona Castorina 110, Rio de Janeiro, 22460-320,
Brasil} \email{henrique@impa.br}
\author{Gil R. Cavalcanti} \address{Mathematical Institute, 24-29 St Giles', Oxford, OX1 3LB, UK.} \email{gilrc@maths.ox.ac.uk}
\author{Marco Gualtieri} \address{Department of Mathematics,
Massachusetts Institute of Technology, 77 Massachussetts avenue,
Cambridge, MA 02139, USA.} \email{mgualt@mit.edu}
\date{}

\begin{abstract}
We develop a theory of reduction for generalized K\"ahler and
hyper-K\"ahler structures which uses the generalized Riemannian
metric in an essential way, and which is not described with
reference solely to a single generalized complex structure.  We show
that our construction specializes to the usual theory of K\"ahler
and hyper-K\"ahler reduction, and it gives a way to view usual
hyper-K\"ahler quotients in terms of generalized K\"ahler reduction.
\end{abstract}

\maketitle

\tableofcontents

\section*{Introduction}

Generalized geometrical structures, such as Dirac structures \cite{Cou90} and
generalized complex or K\"ahler structures \cite{Gu03,Hit03}, are similar to their
classical counterparts, i.e. integrable tangent distributions,
complex or K\"ahler structures, except that they operate not on the
tangent bundle but on the direct sum $TM\oplus T^*M$, or more
generally, on an exact Courant algebroid $E$ which is an extension
of $TM$ by $T^*M$.

In the presence of an action on the underlying manifold by a Lie
group $G$, one may ask whether such generalized geometries
\emph{reduce} to a suitable quotient space.  To address this
question, we developed in~\cite{BCG05} a theory of reduction for
exact Courant algebroids and their associated generalized
geometrical structures, based on the idea of \emph{extending} the
$G$-action on $TM$ to an action on the Courant algebroid $E$.  This
is particularly interesting since $E$ has symmetries beyond the
usual diffeomorphisms; there are also the \emph{B-field}
transformations, well known to physicists.  See \cite{Hu,
LinTol,StXu,Vais} for related work.

In this article, we provide a streamlined approach to the reduction
procedures described in~\cite{BCG05}, focusing on the useful special
case where the extended $G$-action is defined by a
bracket-preserving map $\widetilde{\psi}:\frak{g}\to \Gamma(E)$ and
an equivariant map $\mu:M\to \frak{h}^*$ for $\frak{h}$ a
$\Gg$-module. The map $\mu$ is called a \emph{moment map}, and is
not associated with any geometrical structure but rather to the
extended action itself. We phrase the constructions in this paper
explicitly in terms of the data $(\widetilde{\psi},\frak{h},\mu)$,
showing that specific choices lead to  well-known reductions.

The main construction in this paper is a reduction procedure for
generalized complex structures which goes beyond the natural
generalization of holomorphic or symplectic quotients developed
in~\cite{BCG05}.  For this, we require a $G$-invariant generalized
Riemannian metric $\G$, compatible with the generalized complex
structure $\J$, and such that the Hermitian structure $(\JJ,\G)$ is
compatible with the $G$-action in a suitable sense.  In particular,
we obtain new reduction procedures for generalized K\"ahler and
generalized hyper-K\"ahler structures.

Finally, we prove that this generalized Hermitian reduction does
specialize to the usual K\"ahler and hyper-K\"ahler reductions, and
that it commutes with the forgetful functors taking hyper-K\"ahler
geometry to generalized K\"ahler geometry.

In a sequel to this work~\cite{BCGne}, we apply the theory developed
here to several infinite-dimensional quotients, including the
generalized (hyper-) K\"ahler structure on the moduli space of
instantons on a generalized (hyper-) K\"ahler 4-manifold first
obtained by Hitchin~\cite{Hit05} (see~\cite{MoV} for the
hyper-K\"ahler case), as well as the generalized K\"ahler structures
on certain Lie groups obtained in~\cite{Gu03}.

We would like to thank the organizers of the Poisson 2006
conference, in particular Giuseppe Dito and Yoshiaki Maeda, for
their assistance and hospitality. H.B. thanks CNPq for financial
support, and G.C. thanks EPSRC for financial support.

\section{Generalized geometry and Courant algebroids}

\subsection{Geometry of $TM\oplus T^*M$}\label{subsec:T+T*}

Given a manifold $M$ of dimension $m$, the direct sum $TM\oplus
T^*M$ is equipped with the following canonical structures: a
fiberwise inner product of signature $(m,m)$ given by
\begin{equation}\label{eq:innerprod}
\IP{X+\xi,Y+\eta}:= \eta(X)+\xi(Y), \;\;\; X+\xi, Y+\eta \in
\Gamma(TM\oplus T^*M),
\end{equation}
and the \emph{Courant bracket} \cite{Cou90} on smooth sections of
$TM\oplus T^*M$, defined by
\begin{equation}\label{eq:Cour}
\Cour{X+\xi, Y+\eta}:= [X,Y] + \mathcal{L}_{X}\eta - i_{Y}d\xi.
\end{equation}
A central idea in ``generalized geometry'' \cite{Gu03,Hit03} is
that $TM\oplus T^*M$, together with the operations
\eqref{eq:innerprod} and \eqref{eq:Cour}, should be thought of as
a ``generalized tangent bundle'' to $M$. This leads to a unified
view of various geometrical structures on $M$, as we now briefly
recall.

A \emph{Dirac structure} \cite{Cou90} on $M$ is a Lagrangian
subbundle $L\subset TM\oplus T^*M$ (i.e., $L=L^\perp$, where
$L^\perp$ is the orthogonal complement of $L$ with respect to the
pairing \eqref{eq:innerprod}) whose space of sections is closed
under the Courant bracket:
\begin{equation}\label{eq:integrab}
\Cour{\Gamma(L),\Gamma(L)}\subseteq \Gamma(L).
\end{equation}
Examples of Lagrangian subbundles of $TM\oplus T^*M$ include the
graphs of 2-forms $\omega:TM \to T^*M$ and bivector fields $\pi:T^*M
\to TM$. In these cases the integrability condition
\eqref{eq:integrab} amounts to $d\omega=0$ and $\pi$ being a Poisson
bivector field; other examples can be found in \cite{Cou90}. All
these definitions clearly carry over to the complexified bundle
$(TM\oplus T^*M)\otimes \mathbb{C}$, leading to \emph{complex Dirac
structures}.

We now consider a special class of complex Dirac structures.  A
\emph{generalized complex structure} \cite{Gu03,Hit03} on $M$ is a
complex structure $\J$ on the bundle $TM\oplus T^*M$ which is
orthogonal with respect to \eqref{eq:innerprod}, and such that the
$+i$-eigenbundle $L \subset (TM\oplus T^*M)\otimes \mathbb{C}$
satisfies the integrability condition \eqref{eq:integrab}. The
orthogonality of $\J$ implies that $L=L^\perp$, so $L$ is a complex
Dirac structure also satisfying
\begin{equation}\label{eq:transvers}
L\cap \overline{L}=\{0\}.
\end{equation}
Conversely, any complex Dirac structure on $M$ satisfying
\eqref{eq:transvers} uniquely determines a generalized complex
structure on $M$. Examples of generalized complex structures include
those determined by complex structures $I:TM\to TM$ and symplectic
structures $\omega:TM\to T^*M$, namely
$$
\J_I= \left ( \begin{matrix} -I & 0 \\
                              0& I^*
\end{matrix}\right ),\;\;\;\;
\J_\omega= \left ( \begin{matrix} 0 & -\omega^{-1} \\
                              \omega& 0
\end{matrix}\right ).
$$
One can similarly define generalized versions of Riemannian metrics,
Hermitian structures, as well as K\"ahler and hyper-K\"ahler
structures \cite{Gu03}, see Sections \ref{subsec:genherm} and
\ref{subsec:gkahler}.

Interpreting geometries on $M$ as structures on $TM\oplus T^*M$
introduces two important features. First of all, one can easily
adapt definitions in order to incorporate geometrical structures
which are ``twisted'' by a closed 3-form on $M$: given $H\in
\Omega^3_{cl}(M)$, one replaces the Courant bracket \eqref{eq:Cour}
by the \emph{$H$-twisted Courant bracket} \cite{SeWe01}
\begin{equation}\label{eq:Hcour}
\Cour{X+\xi, Y+\eta}_H= [X,Y] + \mc{L}_X \eta - i_Y d\xi  + i_Y
i_X H,
\end{equation}
and changes the integrability condition \eqref{eq:integrab}
accordingly; this leads to \emph{$H$-twisted Dirac structures},
\emph{$H$-twisted generalized complex structures} and so on. Second,
there is an action of the additive group of closed 2-forms
$\Omega^2_{cl}(M)$ on $TM\oplus T^*M$ preserving the operations
\eqref{eq:innerprod}, \eqref{eq:Cour}, \eqref{eq:Hcour}, given by
$$
X+\xi {\mapsto} X + \xi + i_XB,
$$
for $B \in \Omega^2_{cl}(M)$. As a result, structures on $TM\oplus
T^*M$, such as generalized complex structures, inherit these extra
symmetries, known as \emph{B-field} transformations.

%%%%%%%%%%%%%%%%%%%%%%%%%%%%%%%%%%%%%%%%%%%%%%%%%%%%%%%%%%%%%%%%%%%%%%%%%%%%%%
\subsection{Exact Courant algebroids}\label{exact courant
algebroids}

An axiomatization of the properties of the Courant bracket on
$TM\oplus T^*M$ leads to the general notion of a Courant
algebroid, introduced in \cite{LWX}. This more intrinsic approach
to the Courant bracket will play an important role in the context
of reduction.

A {\it Courant algebroid} over a manifold $M$ is a vector bundle
$\E \to M$ equipped with a fibrewise nondegenerate symmetric
bilinear form $\IP{\cdot,\cdot}$, a bilinear bracket
$\Cour{\cdot,\cdot}$ on the smooth sections $\Gamma(\E)$, and a
bundle map $\pi: \E\to TM$ (called the \emph{anchor}), such that,
for all $e_1,e_2,e_3\in \Gamma(\E)$ and $f\in C^{\infty}(M)$, the
following properties are satisfied:
\begin{itemize}
\item[C1)] $\Cour{e_1,\Cour{e_2,e_3}} = \Cour{\Cour{e_1,e_2},e_3}
+ \Cour{e_2,\Cour{e_1,e_3}}$,

\item[C2)] $\Cour{e_1,fe_2}=f\Cour{e_1,e_2}+
(\mathcal{L}_{\pi(e_1)}f) e_2$,

\item[C3)] $\mathcal{L}_{\pi(e_1)}\IP{e_2,e_3}=
\IP{\Cour{e_1,e_2},e_3} + \IP{e_2, \Cour{e_1, e_3}}$,

\item[C4)] $\pi(\Cour{e_1,e_2})=[\pi(e_1),\pi(e_2)]$,

\item[C5)] $\Cour{e_1,e_1} = \frac{1}{2} \pi^* d\IP{e_1,e_1}$,
\end{itemize}
where in (C5) we identify $E\cong E^*$ via $\IP{\cdot,\cdot}$ in
order to view $\pi^*$ as taking values in $E$. The model example
of a Courant algebroid is $TM\oplus T^*M$, with pairing
\eqref{eq:innerprod}, anchor given by the canonical projection
$TM\oplus T^*M \to TM$, and bracket given by the $H$-twisted
Courant bracket \eqref{eq:Hcour}. It is straightforward to extend
the concepts of Dirac structures, generalized complex structures
etc. to general Courant algebroids.

In the definition of a Courant algebroid, properties C1)--C4)
express natural compatibility conditions between the anchor $\pi$,
the bracket $\Cour{\cdot,\cdot}$ and the pairing $\IP{\cdot,\cdot}$
that will be further discussed in Section \ref{subsec:inner}.
Property C5), on the other hand, prevents the bracket
$\Cour{\cdot,\cdot}$ from being skew-symmetric,
%; instead, we have
%$$
%\Cour{e_1,e_2}=-\Cour{e_2,e_1}+\pi^* d\IP{e_1,e_2}.
%$$
and it implies that $\pi \circ \pi^* =0$, so we have a chain
complex
\begin{equation}\label{exact}
0 \longrightarrow T^*M \stackrel{\pi^*}{\longrightarrow} \E
\stackrel{\pi}{\longrightarrow} TM \longrightarrow 0.
\end{equation}
In this paper we will restrict our attention to \emph{exact
Courant algebroids}, i.e., those for which the sequence
\eqref{exact} is exact. In this case, we always identify $T^*M$
with a subspace of $\E$ via $\pi^*$.

%, so that, for $\xi \in T^*M$ and $e \in \E$,
%$$\IP{\xi, e} = \frac{1}{2}\xi(\pi(e)).$$

Given an exact Courant algebroid $E\to M$, we can always choose a
right splitting $\nabla: TM\to \E$ of \eqref{exact} which is
\emph{isotropic}, i.e., whose image in $\E$ is isotropic with
respect to $\IP{\cdot,\cdot}$. Each such $\nabla$ defines a
``curvature'' 3-form $H\in \Omega^3_{cl}(M)$ by
\begin{equation}\label{eq:curv}
H(X,Y,Z) := \IP{\Cour{\nabla(X),\nabla (Y)},\nabla(Z)},\;\; \mbox{
for } X,Y,Z\in \Gamma(TM).
\end{equation}
Under the vector bundle isomorphism $\nabla+\pi^*:TM\oplus T^*M
{\to} \E$, the Courant algebroid structure on $E$ is identified
with the usual Courant algebroid structure on $TM\oplus T^*M$
defined by the $H$-twisted Courant bracket \eqref{eq:Hcour}.

\begin{remark} Exact Courant algebroids were first studied by P. \v{S}evera,
who classified them by noticing that the choice of a different
isotropic splitting of \eqref{exact} modifies $H$ by an exact
3-form. As a result, the cohomology class $[H]\in H^3(M,\mathbb{R})$
is independent of the splitting and completely determines the exact
Courant algebroid $E$ up to isomorphism. We call $[H]$ the {\it
\v{S}evera class} of $E$.
\end{remark}

\section{Extended actions on Courant algebroids}\label{sec:actions}

In this section we briefly review the notion of \emph{extended
action}, introduced in \cite{BCG05} for the purpose of describing
the reduction of Courant algebroids.

\subsection{Infinitesimal actions}\label{subsec:inner}

The action of a Lie group $G$ on a manifold $M$ may be described
infinitesimally as a Lie algebra homomorphism $\Gg\lra \Gamma(TM)$.
 The definition of an extended action on a Courant algebroid $E$ is
analogous, with $E$ playing the role of $TM$.

Recall that an \emph{infinitesimal automorphism} of  a vector bundle
$E$ is a pair $(F,X)$, where $X\in \Gamma(TM)$ and $F:\Gamma(E)\to
\Gamma(E)$ satisfies
\begin{equation}\label{eq:leibniz}
F(f e)=fF(e)+(\mathcal{L}_Xf)e, \;\;\; e\in \Gamma(E), f \in
C^\infty(M).
\end{equation}
Infinitesimal bundle automorphisms form a Lie algebra with respect
to the bracket
$$
[(F_1,X_1),(F_2,X_2)]:=(F_1F_2-F_2F_1,[X_1,X_2]).
$$
If $E$ is a Courant algebroid, then its Lie algebra of symmetries,
denoted by $\sym(E)$, consists of infinitesimal bundle automorphism
$(F,X)$ which preserve the bracket $\Cour{\cdot,\cdot}$, the pairing
$\IP{\cdot,\cdot}$, and the anchor $\pi:E\to TM$:
\begin{align*}
F(\Cour{e_1,e_2})&=\Cour{F(e_1),e_2} + \Cour{e_1,F(e_2)},\\
\mathcal{L}_X\IP{e_1,e_2}&=\IP{F(e_1),e_2}+\IP{e_1,F(e_2)},\\
\pi\circ F&= \mathcal{L}_X\circ \pi,
\end{align*}
where $e_1,e_2\in \Gamma(E)$. Given a section $e\in \Gamma(E)$, we
observe that axioms C1)--C4) in the definition of a Courant
algebroid imply that the pair $(F=\Cour{e,\cdot},X=\pi(e))$ is in
$\sym(E)$. As a result, we obtain a map:
\begin{equation}\label{eq:Ad2}
\Ad:\Gamma(E)\lra \sym(E),\;\; e \mapsto (\Cour{e,\cdot},\pi(e)).
\end{equation}
The elements of $\sym(E)$ in the image of $\Ad$ are called
\emph{inner symmetries} of $E$. Note that the map \eqref{eq:Ad2}
extends the usual identification of vector fields $X\in \Gamma(TM)$
with infinitesimal symmetries of the Lie bracket on $TM$:
\begin{equation}\label{eq:Ad1}
\Gamma(TM)\longrightarrow \sym(TM),\;\; X\mapsto ([X,\cdot],X).
\end{equation}
It is important to note, however, that although \eqref{eq:Ad1} is
an isomorphism, the map \eqref{eq:Ad2} is neither injective nor
surjective in general.

Given a Lie algebra $\frak{g}$, an equivariant structure on $E$
preserving its Courant algebroid structure is defined
infinitesimally by a Lie algebra homomorphism $\frak{g}\to \sym(E)$.
The particular situation that will concern us in this paper is that
of $\frak{g}$-actions by inner symmetries, i.e. compositions
\begin{equation}\label{eq:factor}
\frak{g}\stackrel{\Psi}{\lra} \Gamma(E)\stackrel{\Ad}{\lra}\sym(E).
\end{equation}
Observe that, since the Courant bracket on $\Gamma(E)$ is \emph{not}
a Lie bracket, it is natural to replace the Lie algebra $\frak{g}$
in~\eqref{eq:factor} by a more general structure with
``Courant-type'' bracket. We call these \emph{Courant
algebras}~\cite{BCG05} and describe them below.

%%%%%%%%%%%%%%%%%%%%%%%%%%%%%%%%%%%%%%%%%%%%%%%%%%%%%%%%%%%%%%%%%%%%%
\subsection{Courant algebras}\label{subsec:couralg}

A \emph{Courant algebra} over a Lie algebra $\Gg$ is a Leibniz
algebra $\Aa$ \cite{Lod93} together with a bracket-preserving map
$\pi:\Aa\lra\Gg$. In other words, $\Aa$ is a vector space endowed
with a bilinear bracket $\Cour{\cdot,\cdot}:\Aa\times\Aa\lra \Aa$
such that
\begin{equation}\label{eq:c1}
\Cour{\ga_1,\Cour{\ga_2,\ga_3}}=\Cour{\Cour{\ga_1,\ga_2},\ga_3}+\Cour{\ga_2,\Cour{\ga_1,\ga_3}},
\;\;\; \ga_1, \ga_2, \ga_3\in\Aa,
\end{equation}
and $\pi(\Cour{\ga_1,\ga_2})=[\pi(\ga_1),\pi(\ga_2)]$ for all
$\ga_1,\ga_2 \in \Aa$. A  Courant algebra is \emph{exact} if $\pi$
is surjective and $\Cour{\gh_1,\gh_2}=0$ for all
$\gh_i\in\ker(\pi)$. In this paper we only consider exact Courant
algebras. Morphisms of Courant algebras are defined in a natural
way.

It is clear that if $E$ is a Courant algebroid, then $\Gamma(\E)$ is
a Courant algebra over the Lie algebra of vector fields, and is
exact if and only if $\E$ is.

If $\Aa\into \Gg$ is an exact Courant algebra, then $\Hh =
\ker(\pi)$ automatically acquires a $\Gg$-module structure: the
action of $\gg\in \Gg$ on $\gh \in \Hh$ is
$$
\gg \cdot \gh := \Cour{\tilde{\gg},\gh},
$$
where $\tilde{\gg}$ is any element of $\Aa$ such that
$\pi(\tilde{\gg}) = \gg$. Since the bracket vanishes on
$\frak{h}$, it easily follows that this is a well defined map
$\Gg\times \Hh \into \Hh$ defining a $\Gg$-action.

The following example of an exact Courant algebra will be central
in this paper.

\begin{example}
[Hemisemidirect product \cite{BCG05,KiWe01}]\label{ex:hemisemi} If
$\frak{g}$ is a Lie algebra and $\frak{h}$ is a $\frak{g}$-module,
then we can endow $\frak{a} :=\frak{g}\oplus \frak{h}$ with the
structure of an exact Courant algebra by taking $\pi:\frak{g}
\oplus \frak{h} \into \frak{g}$ to be the  natural projection and
defining
\begin{equation}\label{eq:hemisemi}
\Cour{(\gg_1,\gh_1),(\gg_2,\gh_2)} := ([\gg_1,\gg_2], \gg_1 \cdot
\gh_2).
\end{equation}
It is clear that $\pi:\Aa \into \Gg$ is surjective and preserves
brackets, and that $\Cour{\Hh,\Hh} = 0$. Finally, condition
\eqref{eq:c1} is a consequence of the Jacobi identity for $\Gg$
and the fact that $\Hh$ is a $\Gg$-module.  This Courant algebra
first appeared in \cite{KiWe01}, where it was studied in the
context of Leibniz algebras.
\end{example}

\subsection{Extended actions}\label{subsec:extend}

Let $E$ be an exact Courant algebroid over $M$. We will now show
how one can produce an infinitesimal action $\frak{g}\to \sym(E)$
starting from an exact Courant algebra $\frak{h}\to
\frak{a}\stackrel{\pi}{\to}\frak{g}$, and a Courant algebra
morphism
$$
\xymatrix{\frak{a} \ar[d]^{\Psi} \ar[r]^{} & \frak{g} \ar[d]^{\psi} \\
\Gamma(E) \ar[r]^{} & \Gamma(TM).}
$$
We will denote a Courant algebra morphism simply by $\Psi:\Aa\to
\Gamma(E)$, keeping in mind that it always projects to an action
on $M$, denoted by $\psi:\frak{g}\to \Gamma(TM)$. It also follows
from the definitions that $\Psi(\frak{h})\subseteq \Omega^1(M)$.

Composing the map \eqref{eq:Ad2} with $\Psi$, we get a map
\begin{equation}\label{eq:map}
\ad\circ \Psi: \frak{a}\lra \sym(E).
\end{equation}
It is important to note that, unlike the map \eqref{eq:Ad1}, the
map $\Ad: \Gamma(E)\to \sym(E)$ has a nontrivial kernel: by
identifying $E$ with $TM\oplus T^*M$ through the choice of an
isotropic splitting, one can directly check that
$\ker(\Ad)=\Omega^1_{cl}(M)$. Hence if $\Psi$ maps $\frak{h}$ into
closed 1-forms, $\Psi(\frak{h})\subseteq \Omega^1_{cl}(M)$, the
map \eqref{eq:map} factors through
$\frak{a}\stackrel{\pi}{\to}\frak{g}$. As a result, there is an
induced infinitesimal $\frak{g}$-action
\begin{equation}\label{eq:actionE}
\frak{g}\lra \sym(E),
\end{equation}
as desired.

We are then led to the following definitions. Let $G$ be a connected
Lie group with Lie algebra $\frak{g}$. An \emph{extended
$\frak{g}$-action} on an exact Courant algebroid $E$ over $M$ is a
Courant algebra morphism $\Psi:\frak{a}\to \Gamma(E)$ from an exact
Courant algebra $\frak{h}\to\frak{a}\stackrel{\pi}{\to}\frak{g}$
into $\Gamma(E)$ so that $\Psi(\frak{h})\subseteq\Omega^1_{cl}(M)$.
We call it an \emph{extended $G$-action} if the induced
$\frak{g}$-action \eqref{eq:actionE} on $E$ integrates to a
$G$-action. In particular, an extended $G$-action on $E$ makes it
into an equivariant $G$-bundle, with $G$ acting by Courant algebroid
automorphisms.

\begin{remark}
Suppose that $\Psi:\frak{a}\to \Gamma(E)$ is an extended
$\frak{g}$-action for which the projected action $\psi:\frak{g}\to
\Gamma(TM)$ integrates to a global $G$-action on $M$. In this
case, a sufficient condition ensuring that this data defines an
extended $G$-action on $E$ (and not only an extended action of a
cover of $G$) is the existence of a $\frak{g}$-invariant isotropic
splitting for ${E}$ (these always exist, e.g., if $G$ is compact),
see \cite[Sec.~2]{BCG05}. Indeed, such a splitting gives an
identification of $E$ with the Courant algebroid $(TM \oplus
T^*M,\IP{\cdot,\cdot},\Cour{\cdot,\cdot}_H)$ satisfying
\begin{equation}\label{H condition}
i_{X_{\ga}}H = d\xi_{\ga},\qquad \mbox{for all } \ga \in \frak{a},
\end{equation}
where $\map(\ga) = X_{\ga} + \xi_{\ga}$. Since $M$ is a
$G$-manifold, $TM\oplus T^*M$ is naturally a $G$-equivariant
bundle, and condition \eqref{H condition} exactly says that this
canonical $G$-action on $TM\oplus T^*M$ coincides infinitesimally
with the one induced by $\Psi$ under the identification $E\cong
TM\oplus T^*M$.
\end{remark}

Any Lie algebra $\Gg$ can be thought of as a Courant algebra over
itself in a trivial way, with the projection $\pi:\frak{g}\to
\frak{g}$ given by the identity. An extended action of this
Courant algebra on an exact Courant algebroid $E$,
$$
\xymatrix{ \Gg \ar[r]^{\Id}\ar[d]^{{\widetilde{\psi}}} & \Gg\ar[d]^{\psi}\\
\Gamma(\E) \ar[r]^{\pi} & \Gamma(TM),}
$$
is called a \emph{lifted (or trivially extended) action} on $E$,
and denoted by $\widetilde{\psi}:\frak{g}\to \Gamma(E)$. A trivial
example is of course when $E=TM\oplus T^*M$, with $H=0$,
and $\widetilde{\psi}=\psi$ is an ordinary
action.

In this paper, we will be only concerned with lifted actions for
which the image of $\widetilde{\psi}:\frak{g} \to E$ is
\emph{isotropic}, i.e., the pairing $\IP{\cdot,\cdot}$ vanishes on
elements $\widetilde{\psi}(u)$.

\begin{example}[Lifted actions]\label{ex:lifted}
Take $E$ to be the Courant algebroid $TM\oplus T^*M$, with $H=0$,
let $\psi:\frak{g}\to \Gamma(TM)$ be an action on $M$, and
$\nu:M\to \frak{g}^*$ be an equivariant map. Then
$$
\widetilde{\psi}(\gg):= \psi(\gg)+d\IP{\nu,u},\;\;\;
\gg\in\frak{g},
$$
is a lifted action on $E$. This observation is a special case of
the following general fact, proven in \cite[Sec.~2]{BCG05}: if we
start with an action $\psi$ on $M$, and assume that $E$ admits an
invariant splitting, with associated 3-form curvature $H$, then
the problem of finding an isotropic lifted action
$\widetilde{\psi}$ extending $\psi$ is equivalent to finding a
closed equivariant extension of $H$ in the Cartan model.
\end{example}

%%%%%%%%%%%%%%%%%%%%%%%%%%%%%%%%%%%%%%%%%%%%%%%%%%%%%%%%%%%%%%%%%%%%%%%%%%%%%%%%%%%%%%%%%%%%%%%%%
\subsection{Moment maps for extended actions}\label{subsec:moment}

Let $\frak{h} \to \frak{a}\to \frak{g}$ be an exact Courant
algebra, and let $\map:\Aa  \to \Gamma(\E)$ be an extended
$\frak{g}$-action on an exact Courant algebroid $E$ over $M$.
Recall that this implies that $\map(\frak{a}) \subset
\Omega_{cl}^1(M)$, and that $\frak{h}$ is a $\frak{g}$-module. Let
us equip $\frak{h}^*$  with the dual $\frak{g}$-action. A
\emph{moment map} for the extended action $\Psi$ is a
$\frak{g}$-equivariant map $\mu: M \into \frak{h}^*$ such that,
for each $\gh \in \Hh$,
$$
\map(\gh) = d\IP{\mu,\gh}.
$$

We now describe how to use equivariant maps to further extend
lifted actions:

\begin{prop}\label{prop:extend}
Let $\widetilde{\psi}:\frak{g}\to \Gamma(E)$ be an isotropic
lifted $\frak{g}$-action on an exact Courant algebroid $E$,
$\frak{h}$ be a $\frak{g}$-module, and $\mu:M\to \frak{h}^*$ be an
equivariant map. Then the map $\Psi:\Gg\oplus \Hh \to \Gamma(E)$,
\begin{equation}\label{eq:extendmoment}
\map(\gg,\gh) = \widetilde{\psi}(\gg) + d\IP{\mu,\gh},
\end{equation}
defines an extended $\frak{g}$-action of the hemisemidirect
product $\Aa = \Gg \oplus \Hh$ on $E$ with moment map $\mu$.
Moreover, the image $\Psi(\frak{a})\subseteq E$ is isotropic over
$\mu^{-1}(0)$.
\end{prop}

\begin{proof}
Let $\Psi$ be defined as in \eqref{eq:extendmoment}. Then, using
that $\Cour{\xi,\cdot}=0$ if $\xi\in \Omega_{cl}^1(M)$, we get
\begin{align*}
\Cour{\map(\gg_1,\gh_1), \map(\gg_2,\gh_2) }&=  \Cour{\widetilde{\psi}(\gg_1),\widetilde{\psi}(\gg_2)}
+\Cour{\widetilde{\psi}(\gg_1),d\IP{\mu,\gh_2}}\\
&=\widetilde{\psi}([\gg_1,\gg_2]) + \mathcal{L}_{\psi(\gg_1)}d\IP{\mu,\gh_2}  \\
&= \map([\gg_1,\gg_2],\gg_1\cdot \gh_2),
\end{align*}
where for the last equality we used the equivariance of $\mu$.
Comparing with \eqref{eq:hemisemi}, we conclude that $\Psi$
preserves brackets. So it defines a Courant algebra morphism. It
is also clear that $\Psi(\frak{h})\subseteq \Omega^1_{cl}(M)$,
hence $\Psi$ is an extended $\frak{g}$-action.

Let us now consider the pairing
$\IP{\Psi(\gg_1,\gh_1),\Psi(\gg_2,\gh_2)}$ over $\mu^{-1}(0)$.
Using that $\widetilde{\psi}$ is isotropic, i.e.,
$\IP{\widetilde{\psi}(\gg_1),\widetilde{\psi}(\gg_2)}=0$, we get
\begin{align*}
\IP{\Psi(\gg_1,\gh_1),\Psi(\gg_2,\gh_2)}&=
\IP{\widetilde{\psi}(\gg_1),d\IP{\mu,\gh_2}} + \IP{ \widetilde{\psi}(\gg_2),d\IP{\mu,\gh_1}}\\
&=\mathcal{L}_{\psi(\gg_1)}\IP{\mu,\gh_2}+
\mathcal{L}_{\psi(\gg_2)}\IP{\mu,\gh_1}\\
&=\IP{\mu,\gg_1\cdot\gh_2} + \IP{\mu,\gg_2\cdot\gh_1},
\end{align*}
which clearly vanishes on points $x\in M$ where $\mu(x)=0$.
\end{proof}

\begin{remark}
Note that if $\widetilde{\psi}$ is a lifted $G$-action for a Lie
group $G$, then $\Psi$ is an extended $G$-action.
\end{remark}

%%%%%%%%%%%%%%%%%%%%%%%%%%%%%%%%%%%%%%%%%%%%%%%%%%%%%%%%%%%%%%%%%%%%%%%%%%%%%%%%%
\section{Reduction of exact Courant algebroids}

In this section, we review the reduction procedure for exact Courant
algebroids introduced in~\cite{BCG05}, giving special attention to
the case described in Prop.~\ref{prop:extend}, i.e. a lifted action
extended by an equivariant map.

\subsection{The general procedure} \label{subsec:courantred}

Let $G$ be a connected Lie group, $\frak{h}\to \frak{a}\to \frak{g}$
be an exact Courant algebra and $\map:\frak{a}\to \Gamma(E)$ be an
extended $G$-action on an exact Courant algebroid $E$. We consider
the distribution
$$
K=\Psi(\frak{a})\subset E
$$
given by the image of the bundle map $\frak{a}\times M \to E$
associated with $\map$.

Our goal is to use the extended action $\Psi$ to produce a
``smaller'' Courant algebroid. To this end, let us suppose that $P
\hookrightarrow M$ is a $G$-invariant submanifold of $M$ such that
\begin{enumerate}
\item $K|_P$ is isotropic;

\item $TP=\pi(K^\perp)$.
\end{enumerate}
One can directly check that $\pi(K^\perp) = \Ann(\map(\frak{h}))$,
so the last condition can be re-written as
\begin{equation}\label{eq:P}
TP=\Ann(\map(\frak{h}))|_P.
\end{equation}

Let us assume the following regularity conditions: $K|_P$ is a
vector bundle over $P$ and the $G$-action on $P$ is free and
proper. As observed in \cite[Sec.~3.1]{BCG05}, it follows from the
fact that $\Psi$ is a Courant algebra morphism that the vector
bundles $K|_P$ and $K^\perp|_P$ are $G$-equivariant subbundles of
$E|_P$. It is then proven in \cite[Thm.~3.3]{BCG05} that the
vector bundle
\begin{equation}\label{eq:Ered}
E^{red}:= \left. \frac{K^\perp|_P}{K|_P}\right / G
\end{equation}
defines an exact Courant algebroid over the manifold $P/G$, called
the \emph{reduced Courant algebroid}.

\begin{remarks}
\item[a)] {As shown in \cite{BCG05}, the assumption that $K|_P$ is
isotropic is not necessary, and this is relevant for some
examples. Note also that $K|_P$ is a vector bundle if and only if
the distribution $\Psi(\frak{h})\subset E$ has constant rank over
$P$.}

\item[b)] The reduced Courant bracket on $E^{red}$ is obtained
canonically, by noticing that the restriction of the bracket on
$E$ to the space $\Gamma(K^\perp)^G$ of $G$-invariant sections of
$K^\perp$ is well-defined modulo $\Gamma(K)^G$, and $\Gamma(K)^G$
is an ideal in $\Gamma(K^\perp)^G$.

\item[c)] Although the reduced Courant algebroid $E^{red}$ is
exact, it may not have a canonical splitting. Nevertheless, one can
still describe the \v{S}evera class of $E^{red}$, see
\cite[Sec.~3.2]{BCG05}.
\end{remarks}

%%%%%%%%%%%%%%%%%%%%%%%%%%%%%%%%%%%%%%%%%%%%%%%%%%%%%%%%%%%%%%%%%%%%%%%%%%%
\subsection{A special case of moment map reduction}\label{subsec:momentreduc}

Let $E$ be an exact Courant algebroid over $M$. We will now
specialize the reduction procedure for actions arising as in
Prop.~\ref{prop:extend}.

\begin{definition}\label{def:data}
Let  $\widetilde{\psi}:\frak{g}\to \Gamma(E)$ be an isotropic lifted
$G$-action on $E$ and let $\mu:M\to \frak{h}^*$ be an equivariant
map, for $\mathfrak{h}$ some $\Gg$-module. Assume that $0$ is a
regular value of $\mu$, and that the $G$-action on $\mu^{-1}(0)$ is
free and proper.  We refer to the triple $(\widetilde\psi,\Hh,\mu)$
as \emph{reduction data}.
\end{definition}

It follows from Prop.~\ref{prop:extend} that reduction data define
an extended $G$-action $\Psi$ of the hemisemidirect product
$\frak{g}\oplus \frak{h}$ on $E$, with image
\begin{equation}\label{eq:K}
K= \{\widetilde{\psi}(\gg)+d\IP{\mu,\gh},\; \gg\in \frak{g},\gh\in
\frak{h}\} \subseteq E.
\end{equation}

\begin{prop}\label{prop:courantred}
Let $E$ be an exact Courant algebroid over $M$, and let
$(\widetilde{\psi},\frak{h},\mu)$ be reduction data. Then
$K|_{\mu^{-1}(0)}$ is an equivariant $G$-bundle over $\mu^{-1}(0)$,
and the quotient vector bundle
$$
\left. E^{red}:=
\frac{K^\perp|_{\mu^{-1}(0)}}{K|_{\mu^{-1}(0)}}\right / G
$$
defines an exact Courant algebroid over $\mu^{-1}(0)/G$.
\end{prop}

\begin{proof}
Note that
$\Ann(\Psi(\frak{h}))|_{\mu^{-1}(0)}=\Ann(d\mu)|_{\mu^{-1}(0)}=T(\mu^{-1}(0))$,
so \eqref{eq:P} holds for $P=\mu^{-1}(0)$. The fact that
$K|_{\mu^{-1}(0)}$ is isotropic follows from
Prop.~\ref{prop:extend}. Using that the $G$-action on $\mu^{-1}(0)$
is free and that $\Psi(\frak{h})|_{\mu^{-1}(0)}$ is a bundle, one
concludes that $K|_{\mu^{-1}(0)}$ is a vector bundle. The result now
follows from the construction outlined in Section
\ref{subsec:courantred}.
\end{proof}

\begin{example}\label{ex:courred}
Let $(TM\oplus T^*M,\IP{\cdot,\cdot},\Cour{\cdot,\cdot}_H)$ be the
Courant algebroid $E$. If the closed 3-form $H$ is basic with
respect to a $G$-action on $M$, then we may regard this $G$-action
as a lifted action on $E$, i.e., $\widetilde{\psi}=\psi:\frak{g}\to
\Gamma(TM)\subset \Gamma(E)$. In this case, for any choice of
$\Gg$-module $\frak{h}$ and equivariant map $\mu:M\to \frak{h}^*$,
the reduced Courant algebroid is naturally split, as follows.  Since
$K=\psi(\frak{g})\oplus d\IP{\mu,\frak{h}}$ and $K^\perp= T
\mu^{-1}(0) \oplus \Ann(\psi(\frak{g}))$, we have
$$
E^{red}= \left. \frac{T \mu^{-1}(0)}{\psi(\frak{g})} \oplus
\frac{\Ann(\psi(\frak{g}))}{d\IP{\mu,\frak{h}}} \right / G =
TM^{red}\oplus T^*M^{red},
$$
where $M^{red}=\mu^{-1}(0)/G$. The curvature 3-form of $E^{red}$
with respect to this natural splitting is just the pushdown of the
(basic) 3-form $\iota^*H$, where $\iota:\mu^{-1}(0)\hookrightarrow
M$ is the inclusion.

This example has two simple extreme cases: if we take
$\frak{h}=\{0\}$, then $E^{red}=T(M/G)\oplus T^*(M/G)$, whereas if
we pick $\frak{g}=\{0\}$, then $E^{red}=T\mu^{-1}(0)\oplus
T^*\mu^{-1}(0)$.
\end{example}

\begin{remarks}
\item[a)] In general, the reduced Courant algebroid $E^{red}$ will not
be canonically split, but one can construct splittings by choosing
connections for the $G$-bundle $\mu^{-1}(0)\into \mu^{-1}(0)/G$,
see \cite[Sec.~3]{BCG05} for a detailed discussion.

\item[b)] A description of the \v{S}evera class of a reduced Courant
algebroid is presented in \cite[Sec.~3.2]{BCG05}, including
explicit examples where a trivial \v{S}evera class is reduced to
a non-trivial one.
\end{remarks}

%%%%%%%%%%%%%%%%%%%%%%%%%%%%%%%%%%%%%%%%%%%%%%%%%%%%%%%%%%%%%%%%%%%%%%%%%
\section{Reduction of geometrical structures}

In this section we explain how various geometrical structures can be
transported to the reduced Courant algebroid in the presence of an
extended action. For simplicity, we restrict our attention to the
special case of extended actions determined by the reduction data
$(\widetilde{\psi},\frak{h},\mu)$, as in
section~\ref{subsec:momentreduc}.

\subsection{Reduction of Dirac and generalized complex
structures}\label{subsec:redgcs}

As in Section \ref{subsec:momentreduc}, $E$ is an exact Courant
algebroid over $M$, $(\widetilde{\psi},\frak{h},\mu)$ defines
reduction data, $K$ is given by \eqref{eq:K}, and the reduced
Courant algebroid $E^{red}$ over $M^{red}=\mu^{-1}(0)/G$ is given in
Prop.~\ref{prop:courantred}.

Suppose that $L\subset E$ is a $G$-invariant Dirac structure,
defining a $G$-equivariant subbundle of $E$. (Infinitesimally,
this means that $\Cour{\widetilde{\psi}(\gg),\Gamma(L)}\subseteq
\Gamma(L)$.) Consider the following distribution of $E^{red}$:
\begin{equation}\label{eq:Lred}
L^{red} :=\left. \frac{(L \cap \Kperp +
K)|_{\mu^{-1}(0)}}{K|_{\mu^{-1}(0)}}\right /G.
\end{equation}
One can directly verify that $L^{red\perp} = L^{red}$. So, at each
point, \eqref{eq:Lred} defines a Lagrangian subspace of $E^{red}$.
However, the distribution $L^{red}$ may not form a smooth vector
bundle over $M^{red}$. (A sufficient, but not necessary, condition
is that $L \cap K |_{\mu^{-1}(0)}$ has constant rank.) If it does,
then it is shown in \cite[Sec.~4.1]{BCG05} that $L^{red}$ is
automatically integrable with respect to the reduced Courant
bracket, and hence defines a Dirac structure in $E^{red}$. The same
construction holds for complex Dirac structures if we replace $K$ by
its complexification $K\otimes \mathbb{C}\subseteq E\otimes
\mathbb{C}$ in \eqref{eq:Lred}, yielding a reduced complex Dirac
structure in $E^{red}\otimes \mathbb{C}$.

If now $\J$ is a $G$-invariant generalized complex structure on $E$,
then its $+i$-eigenbundle $L$ is a $G$-invariant complex Dirac
structure in $E\otimes \mathbb{C}$. So one may attempt to transport
$\J$ to the reduced Courant algebroid $E^{red}$ by reducing this
Dirac structure as in \eqref{eq:Lred}. Assuming that $L^{red}$ is a
smooth bundle over $M^{red}$, then it defines a reduced generalized
complex structure $\J^{red}$ in $E^{red}$ if and only if it
satisfies the extra condition
\begin{equation}\label{eq:redtrans}
 L^{red} \cap \overline{L^{red}} =
\{0\}.
\end{equation}
As observed in \cite[Sec.~5.1]{BCG05}, this last condition can be
equivalently expressed in terms of the operator $\J$ as follows:
over each point in $\mu^{-1}(0)$, we must have
\begin{equation}\label{eq:reduction condition}
\J K \cap \Kperp \subset K.
\end{equation}
A particularly simple condition implying both the smoothness of
$L^{red}$ and condition \eqref{eq:reduction condition} is
\begin{equation}\label{eq:JKK}
\J K=K \;\;\; \mbox{ over } \mu^{-1}(0),
\end{equation}
as discussed in \cite[Sec.~5]{BCG05}. We summarize the discussion by
citing the following result:
\begin{theo}[\cite{BCG05}, Thm.~5.2]
Let $\J$ be a generalized complex structure on the exact Courant
algebroid $E$ and let $(\widetilde{\psi},\frak{h},\mu)$ be reduction
data. If $\J$ is $G$-invariant and satisfies \eqref{eq:JKK}, then
$L^{red}$ defines a reduced \gcs\ $\J^{red}$ on $E^{red}$.
\end{theo}
We now illustrate this reduction procedure with two simple examples,
in which the Courant algebroid $E$ is taken to be $TM\oplus T^*M$,
with $H=0$.

\begin{example}[Hamiltonian reduction]\label{ex:hamreduc}
Consider a hamiltonian $G$-manifold $(M,\omega)$, with action
$\psi:\frak{g}\to \Gamma(TM)$ and moment map $\mu:M\to \frak{g}^*$.
We can describe Hamiltonian reduction as generalized reduction as
follows: the triple $(\psi,\frak{g},\mu)$ defines reduction data,
and, according to Example~\ref{ex:courred}, the reduced Courant
algebroid is $E^{red}=TM^{red}\oplus T^*M^{red}$, where
$M^{red}=\mu^{-1}(0)/G$. Viewing $\omega$ as a generalized complex
structure $\J_{\omega}$, it follows from the \emph{moment map
condition}
\begin{equation}\label{eq:mmapcond}
i_{\psi(\gg)}\omega=d\IP{\mu,\gg},\;\;\;
\gg\in \frak{g}
\end{equation}
that the reduction data and the \gcs\ are related by
\begin{equation}\label{eq:comp}
\J_{\omega} d\IP{\mu,\gg}=\psi(\gg),
\end{equation}
and this immediately implies condition \eqref{eq:JKK}. So we can
carry out generalized reduction to obtain
$\J_{\omega}^{red}=\J_{\omega_{red}}$, where $\omega_{red}$ is the
reduced symplectic form on $M^{red}$ obtained by Marsden-Weinstein
reduction.
\end{example}

\begin{remark}
The compatibility \eqref{eq:comp} of the previous example can be
generalized as follows: instead of ordinary actions $\psi$, one can
consider more general lifted actions, for example those of the form
$\widetilde{\psi}(\gg)=\psi(\gg)+d\IP{\nu,\gg}$, where $\nu:M\to
\frak{g}^*$ is equivariant (see Example \ref{ex:lifted}); one can
also consider more general $\Gg$-modules $\frak{h}$, and then impose
the condition
$$
\J d\IP{\mu,\gg}=\widetilde{\psi}(\gg)=\psi(\gg)+d\IP{\nu,\gg},
$$
which directly implies \eqref{eq:JKK}. If $\frak{h}=\frak{g}$, these
are the actions studied in \cite{Hu,LinTol} (where the map
$\mu+i\nu$ is called the moment map). For the examples of interest
in this paper, we will need $\frak{g}\neq \frak{h}$ and a weaker
version of \eqref{eq:JKK} to hold (c.f. Example~\ref{ex:hk}).
\end{remark}

\begin{example}[Holomorphic quotients]\label{ex:holreduc}
Suppose that a complex group $(G,I_G)$ acts holomorphically on a complex
manifold $(M,I)$. We now consider the reduction data $(\psi,
\frak{h}=\{0\}, \mu=0)$. In this case, the reduced Courant algebroid
is $T(M/G)\oplus T^*(M/G)$. On the other hand, condition
\eqref{eq:JKK} is a direct consequence of the fact that the action
is holomorphic: $\psi(I_Gu)=I\psi(u)$, $u \in \frak{g}$.
Carrying out generalized reduction for $\J_I$, we obtain $\J_I^{red}
= \J_{I^{red}}$, where $I^{red}$ is the quotient complex structure
on $M/G$.
\end{example}

The framework of reduction described in this section is general
enough to include more exotic examples (see e.g.
\cite[Sec.~5.2]{BCG05}), such as the case of a lifted action
preserving a symplectic structure whose reduction is a complex
structure.

%%%%%%%%%%%%%%%%%%%%%%%%%%%%%%%%%%%%%%%%%%%%%%%%%%%%%%%%%%%%%%%%%%%%%%
\subsection{Generalized Hermitian reduction}\label{subsec:genherm}

We now describe a situation in which a generalized complex structure
may have a natural reduction even when condition~\eqref{eq:JKK}
fails to hold. While we state our results for extended actions
defined by reduction data $(\widetilde{\psi},\frak{h},\mu)$, this is
not essential; more general extended actions may be treated in the
same way.

 A {\it generalized (Riemannian) metric} on a Courant algebroid $\E$ is an
orthogonal, self-adjoint bundle automorphism $\G:\E \into \E$
satisfying
$$
\IP{\G e,e} > 0,\;\;\; \forall e\neq 0.
$$
A generalized metric is compatible with a \gcs\ $\J$ if they
commute. The pair $(\J,\G)$ is then called a {\it generalized
Hermitian structure}.

If $K\subset E$ and $\G$ is a generalized metric on $E$, then we
define
\begin{equation}
K^{\G} := \G \Kperp \cap \Kperp,
\end{equation}
which is the $\G$-orthogonal complement of $K$ in $\Kperp$.

\begin{theo}[Generalized Hermitian reduction]\label{theo:gcss II}
Let $E$ be an exact Courant algebroid over $M$, with reduction data
$(\widetilde{\psi},\frak{h},\mu)$. Suppose that $E$ is equipped with
a $G$-invariant generalized Hermitian structure $(\J,\G)$. If
\begin{equation}\label{eq:JKKweak}
\J K^{\G} =K^{\G}\;\; \mbox{ over } \mu^{-1}(0),
\end{equation}
then $\J$ can be reduced to $E^{red}$, and $\G$ induces a compatible
generalized metric on  $E^{red}$.
\end{theo}

\begin{proof}
The proof follows closely the ideas in \cite[Sec.~6.1]{BCG05}. Let
us first notice that condition \eqref{eq:JKKweak} implies
\eqref{eq:reduction condition}: Since $\J$ is orthogonal, we see
that
 $$
 \IP{\J K,K^{\G}}=\IP{K, K^{\G}} = 0,
 $$
 hence $\J K \subset K^{\G \perp}$ and $\J K \cap \Kperp \subset  K^{\G \perp} \cap \Kperp =
 K$. This shows that the Dirac reduction $L^{red}$ of the
 $i$-eigenbundle of $\J$ satisfies \eqref{eq:redtrans}, hence it
 defines a reduced \gcs\ as long as it is a smooth bundle, a fact
 which we now verify.

Using the $\G$-orthogonal decomposition $K^\perp=K\oplus K^{\G}$
over $\mu^{-1}(0)$, we obtain an identification of vector bundles
$$
K^{\G}|_{\mu^{-1}(0)}/G \cong \E^{red}.
$$
Since $K^{\G}$ is $\J$ invariant, this identification induces a
\gacs\ $\J^{red}$ on $E^{red}$, whose $+i$-eigenbundle agrees with
the reduced Dirac structure $L^{red}$. This implies that $L^{red}$
is smooth, and hence integrable, and that $\J^{red}$ is the \gcs\
associated to it.

Since $\G(K^{\G}) = K^{\G}$, we can also transport the generalized
metric $\G$ to a generalized metric $\G^{red}$ on $\E^{red}$, and
since $\G$ and $\J$ commute pointwise, the same holds for their
restrictions to $K^{\G}$. Thus the reduced metric and \gcs\ are
compatible.
\end{proof}

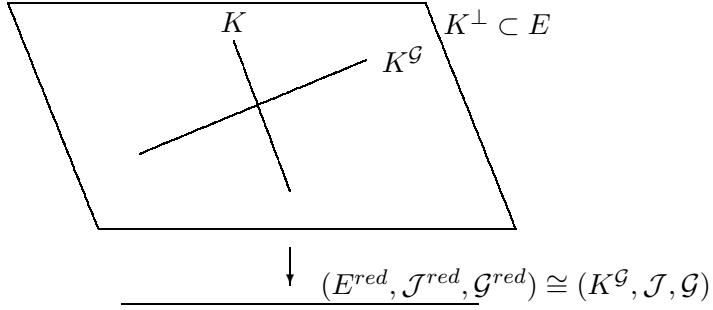
\begin{figure}[h!]
\begin{center}
\unitlength 0.5mm
\begin{picture}(140,90)(0,0)
\linethickness{0.1mm}
\multiput(5,90)(0.12,-0.3){200}{\line(0,-1){0.3}}
\put(5,90){\line(1,0){111}} \put(29,30){\line(1,0){111}}
\multiput(116,90)(0.12,-0.3){200}{\line(0,-1){0.3}}
%\linethickness{0.1mm} \qbezier(50,80)(74.22,70.62)(73.12,53.75)
%\qbezier(73.12,53.75)(72.03,36.88)(70,35) \linethickness{0.2mm}
\multiput(65,80)(0.12,-0.32){125}{\line(0,-1){0.32}}
%\put(35,85){\makebox(0,0)[cc]{``$G$-orbit''}}

\put(65,85){\makebox(0,0)[cc]{$K$}}

\put(110,75){\makebox(0,0)[cc]{$K^{\G}$}}

\put(135,85){\makebox(0,0)[cc]{$K^{\perp}\subset E$}}

\linethickness{0.1mm} \put(80,15){\line(0,1){10}}
\put(80,15){\vector(0,-1){0.12}}
%\put(90,20){\makebox(0,0)[cc]{$p$}}
\put(90,20){\makebox(0,0)[cc]{}}

\put(140,15){\makebox(0,0)[cc]{$(\E^{red},\J^{red},\G^{red}) \cong
(K^{\G},\J,\G)$}}

\linethickness{0.1mm}
\multiput(40,50)(0.29,0.12){208}{\line(1,0){0.29}}
\linethickness{0.1mm} \put(35,10){\line(1,0){95}}
\end{picture}
\caption{The generalized metric and complex structure on the reduced
Courant algebroid are modeled on the $\G$-orthogonal complement of
$K$ inside $\Kperp$.} \label{fig:algebroid viewpoint}
\end{center}
\end{figure}

Intuitively, thinking of $E$ as a generalized tangent bundle to $M$,
the distribution $K$ plays the role of the tangent distribution to
the orbits of the extended action. The construction above can be
interpreted as saying that the orthogonal complement of the
generalized ``$G$-orbit'' is a model for the quotient; see Figure
\ref{fig:algebroid viewpoint}. Despite the clarity of condition
\eqref{eq:JKKweak}, it is often easier to verify its orthogonal
complement:
\begin{equation}\label{eq:easier condition}
\J (K +\G K) = K +\G K.
\end{equation}
For further details concerning the reduction of generalized metrics,
see~\cite{Ca06c}.

%%%%%%%%%%%%%%%%%%%%%%%%%%%%%%%%%%%%%%%%%%%%%%%%%%%%%%%%%%%%%%%%%%%%%%%
\subsection{Generalized K\"ahler and hyper-K\"ahler
reductions}\label{subsec:gkahler}

A {\it \gks} is a \gcs\ $\J$ together with a compatible generalized
Riemannian metric $\G$ such that $ \J' := \J\G$ is also a \gcs.

A {\it \ghk} structure is a triple of \gcss, $\J_1,\J_2,\J_3$, each
of which forms a generalized K\"ahler structure with the same
generalized Riemannian metric $\G$, and such that
\[
\J_1 \J_2 = -\J_2\J_1=\J_3.
\]

\begin{remark}
Observe that given a \gcs\ $\J$ with compatible metric $\G$,
the product $\J \G$ is always a \gacs, i.e., $(\J\G)^2 = -\Id$.
The nontrivial requirement is the integrability of this structure.
Note also that $\J \G$ is compatible with $\G$.
\end{remark}

We may now apply Theorem~\ref{theo:gcss II} to obtain \gk\ and \ghk\
reductions. As before, let $(\widetilde{\psi},\frak{h},\mu)$ be
reduction data for the exact Courant algebroid $E$.

\begin{theo}[Generalized K\"ahler reduction]\label{theo:gk}
Suppose that $(\J,\G)$ is a $G$-invariant \gks\ on $E$. If $\J
K^{\G} =K^{\G}$ over $\mu^{-1}(0)$, then the generalized K\"ahler
structure $(\J,\G)$ reduces to the Courant algebroid $E^{red}$ over
$\mu^{-1}(0)/G$.
\end{theo}
\begin{proof}
Letting $\J' = \J \G$, it is clear that $\J'$ is also
$G$-invariant and that both $(\J,\G)$ and $(\J',\G)$ satisfy the
conditions of  Theorem \ref{theo:gcss II}. Hence they can be
reduced to $E^{red}$, and the reduced metric $\G^{red}$ is
compatible with both reduced generalized complex structures. Since
the reduced structures are identified with the restrictions of
$\J,\J' $ and $\G$ to $K^{\G}$, the fact that $\J' = \J \G$ on
$K^{\G}$ implies that the same holds in $\E^{red}$. Hence
$(\J^{red},\G^{red})$ is a \gks.
\end{proof}

We have the following analogous result for \ghk\ structures:

\begin{theo}[Generalized hyper-K\"ahler reduction]\label{theo:ghk}
Suppose that $(\J_1,\J_2,\J_3,\G)$ is a $G$-invariant \ghks\ on $E$.
If $\J_1 K^{\G}= \J_2 K^{\G} =\J_3 K^{\G} =K^{\G}$ over
$\mu^{-1}(0)$, then the \ghks\ $(\J_1,\J_2,\J_3,\G)$ reduces to a
\ghks\ in $E^{red}$ over $\mu^{-1}(0)/G$.
\end{theo}

%%%%%%%%%%%%%%%%%%%%%%%%%%%%%%%%%%%%%%%%%%%%%%%%%%%%%%%%%%%%%%%%
\subsection{Examples}

We now describe how one recovers the classical K\"ahler and
 hyper-K\"ahler quotients \cite{GuSt82,HKLR87,Kir84} using our methods.
 As a final example we show that the
classical hyper-K\"ahler quotient may be viewed as a non-trivial
example of a \gk\ quotient.

For the examples in this section, the Courant algebroid in question
is $E=TM\oplus T^*M$, with $H=0$, and with reduction data
$(\psi,\frak{h},\mu)$, where the lifted action $\psi$ is an ordinary
action.

\begin{remark}
In a separate paper~\cite{BCGne}, we describe examples of
generalized K\"ahler and hyper-K\"ahler quotients involving
non-vanishing twists $H$ and non-trivial lifted actions, such as the
construction of generalized K\"ahler and hyper-K\"ahler structures
on certain moduli spaces of instantons and Lie groups.
\end{remark}

\begin{example}[K\"ahler reduction]\label{ex:kahler}
Let $(I,\omega)$ be a K\"ahler structure on $M$, preserved by a
Hamiltonian $G$-action $\psi:\frak{g}\to \Gamma(TM)$, with moment
map $\mu:M\to \frak{g}^*$. Let $g=\omega I$ be its associated
K\"ahler metric.

We then view the K\"ahler structure as a generalized K\"ahler
structure $(\J_{\omega},\G)$, where
\begin{equation}\label{eq:genmetric}
\G=\begin{pmatrix}
0 & g^{-1}\\
g & 0
\end{pmatrix},
\end{equation}
and carry out its reduction using the reduction data
$(\psi,\frak{g},\mu)$. The reduced Courant algebroid, as described
in Example~\ref{ex:courred}, is simply $TM^{red}\oplus T^*M^{red}$,
where $M^{red}=\mu^{-1}(0)/G$.

As in Example \ref{ex:hamreduc}, the moment map condition
\eqref{eq:mmapcond} for $\omega$ immediately implies that
$\J_\omega K=K$. Hence
 $K^{\G}=\J_I \J_{\omega}\Kperp  \cap  K^{\perp} = \J_I \Kperp \cap \Kperp$ and therefore
$$
\J_\omega K^\G = \J_{\omega}(\J_I \Kperp \cap \Kperp) = \J_I
\Kperp \cap \Kperp = K^\G.
$$
So we can apply  Theorem \ref{theo:gk} and reduce the K\"ahler
structure as a generalized K\"ahler structure. We now check that
the reduced structure $\J_I^{red}$ and $\J_{\omega}^{red}$ agree
with the generalized structures associated with the usual reduced
K\"ahler structure on $M^{red}$. First, as observed in
Example~\ref{ex:hamreduc}, the generalized reduction of the
symplectic structure is just the usual symplectic reduction. Let
us now discuss the reduction of $\J_I$. We know that, at each
point of $M^{red}$, $\J_I^{red}$ is described by its restriction
to
$$
K^{\G} = \{X +\xi \in TP \oplus T^*P:X \perp \psi(\frak{g}) \mbox{
and } \xi(\psi(\Gg)) =0\},
$$
where $P=\mu^{-1}(0)$. We have the identifications
$$
\{X \in TP :X \perp \psi(\frak{g})\} = TM^{red}\;\; \mbox{ and
}\;\; \{\xi \in T^*P:\xi(\psi(\Gg)) =0\} = T^*M^{red},
$$
and the space on the left is invariant under $\J_I$. It follows that
$\J_I^{red}$ preserves $TM^{red}$, and hence it is of complex type.
One sees from this description that $\J_I^{red}$ agrees with the
\gcs\ associated with the complex structure obtained by the usual
K\"ahler reduction procedure \cite{GuSt82,Kir84}.
\end{example}

\begin{example}[hyper-K\"ahler reduction]\label{ex:hk}
Let $(I_1,I_2,I_3,g)$ be a hyper-K\"ahler structure on $M$ preserved
by a Hamiltonian $G$-action $\psi:\frak{g}\to \Gamma(TM)$, in the
sense that there exist moment maps $\mu_1,\mu_2,\mu_3 \in
C^\infty(M,\frak{g}^*)^G$ satisfying
\begin{equation}\label{eq:hmmapcond}
i_{\psi(\gg)} \omega_j = d\IP{\mu_j,\gg}, \;\; \forall \gg \in \Gg,
\mbox{ and } j = 1,2,3,
\end{equation}
where $\omega_j=gI_j$ are the K\"ahler forms.

We consider the generalized complex structures $\J_j$ associated
with the complex structures $I_j$, $j=1,2,3$, and the generalized
metric $\G$ as in \eqref{eq:genmetric}. Then $(\J_1,\J_2,\J_3,\G)$
defines a generalized hyper-K\"ahler structure, and we consider the
reduction data $(\psi,\frak{h},\mu)$, where now
$$
\frak{h}=\frak{g}\oplus \frak{g}\oplus \frak{g},\;\;\;\mbox{ and
}\;\;\; \mu=(\mu_1,\mu_2,\mu_3):M\to \frak{h}^*.
$$
It follows that
$$
K = \{\psi(\gg) + d\IP{\mu_1,\gh_1} + d\IP{\mu_2,\gh_2}
+d\IP{\mu_3,\gh_3},\; \gg\in \frak{g}, \;\; \gh_1,\gh_2,\gh_3 \in
\frak{h} \}.
$$
In order to apply Theorem \ref{theo:ghk}, we must check that
\begin{equation}\label{eq:checking}
\J_j (K + \G K) = K + \G K
\end{equation}
over $\mu^{-1}(0)$, $j=1,2,3$. Using \eqref{eq:hmmapcond}, we see
that
$$
K  + \G K =\psi(\Gg) +
\omega_1(\psi(\Gg))+\omega_2(\psi(\Gg))+\omega_3(\psi(\Gg)) +
g(\psi(\Gg)) + I_1 \psi(\Gg) + I_2 \psi(\Gg)+ I_3 \psi(\Gg).
$$
It now follows from the relations $I_1 I_2 = I_3$, $I_1 \omega_2
=- \omega_3$ and $I_1 \omega_1 = -g$ that \eqref{eq:checking}
holds. Hence we can reduce the hyper K\"ahler-structure to
$M^{red}=\mu^{-1}(0)/G$ as a generalized hyper-K\"ahler structure.

As in Example~\ref{ex:kahler}, we have the identification
\begin{equation}\label{eq:KG}
K^{\G} =\{X \in TP :X \perp \psi(\frak{g})\}\oplus \{\xi \in
T^*P:\xi(\psi(\Gg)) =0\} = TM^{red} \oplus T^*M^{red},
\end{equation}
for $P=\mu^{-1}(0)$. Since $K^{\G}$ is invariant under $\J_1$,
$\J_2$ and $\J_3$ and these structures are of complex type, it
follows that the space
$$
\{X \in TP :X \perp \psi(\frak{g})\} \cong TM^{red},
$$
is also invariant by these structures. Hence
$\J_1^{red},\J_2^{red}$ and $\J_3^{red}$ are complex structures
and the \ghks\ obtained in the reduced manifold is precisely the
usual hyper-K\"ahler reduction of $M$ from \cite{HKLR87}.
\end{example}

In~\cite{Gu03}, it was shown that a generalized K\"ahler structure
$(\JJ,\G)$ on $E=TM\oplus T^*M$ with $H=0$ determines and is
uniquely determined by a quadruple $(I_+,I_-,g,b)$, where $g$ is a
Riemannian metric on $M$, $I_+$ and $I_-$ are Hermitian complex
structures (hence defining a \emph{bihermitian} structure), and $b$
is a 2-form such that
\begin{equation}\label{eq:gk bihermitian}
d^c_-\omega_-=- d^c_+\omega_+ = db,
\end{equation}
where $\omega_{\pm}=gI_{\pm}$ and $d^c=i(\delbar -
\partial)$ is defined by the appropriate complex structure.

The bihermitian structure is obtained as follows (see~\cite{Gu03}
for details). Since $\G^2=\Id$, we can write $E=C_+\oplus C_-$,
where $C_{\pm}$ is the $\pm 1$-eigenbundle of $\G$. The spaces $C_+$
and $C_-$ intersect $T^*M$ trivially, so the projection $\pi:E\to
TM$ induces identifications $C_\pm \cong TM$. The metric $g$ on $M$
is induced by the restriction of $\IP{\cdot,\cdot}$ to $C_+$,
whereas $J_\pm$ come from the restrictions of $\J$ to $C_{\pm}$.
Conversely, the generalized K\"ahler structure may be written in
terms of $(I_+,I_-,g,b)$ as follows.
\begin{gather*}
\JJ=\frac{1}{2}\left(\begin{matrix}1&\\b&1\end{matrix}\right)
\left(\begin{matrix}I_++ I_- & -(\omega_+^{-1}-\omega_-^{-1}) \\
\omega_+-\omega_-&-(I^*_++I^*_-)\end{matrix}\right)
\left(\begin{matrix}1&\\-b&1\end{matrix}\right),\\
G=\left(\begin{matrix}1&\\b&1\end{matrix}\right)
\left(\begin{matrix}&g^{-1}\\g&\end{matrix}\right)
\left(\begin{matrix}1&\\-b&1\end{matrix}\right).
\end{gather*}

It follows from this bihermitian interpretation of generalized
K\"ahler geometry that any hyper-K\"ahler structure $(g, I_1, I_2,
I_3)$ determines a generalized K\"ahler structure, by simply
choosing $I_+=I_1$ and $I_-=I_2$, for example.  We now show that
this ``forgetful functor'' commutes with reduction, i.e.
intertwines the notion of generalized K\"ahler reduction from
Theorem~\ref{theo:gk} with the usual hyper-K\"ahler quotient
procedure.

\begin{example}[Hyper-K\"ahler reduction versus generalized K\"ahler reduction]\label{ex:gk}
A hyper-K\"ahler structure $(g,I_1,I_2,I_3)$ defines a bihermitian
structure $(g,I_1,I_2)$ satisfying \eqref{eq:gk bihermitian} for
$b=0$, hence it defines a generalized K\"ahler structure $(\G,\J)$
as above.  The \gcs\ $\J$ may be described as a bundle automorphism
of $E=C_+\oplus C_-$ as follows: on $C_+$,
$$
\J(X + g(X)) = I_1X + g(I_1X),
$$
and on $C_-$,
$$
\J(X - g(X)) = I_2 X -g(I_2X).
$$
Given a Hamiltonian $G$-action $\psi:\frak{g} \into \Gamma(TM)$
preserving the hyper-K\"ahler structure, and with moment maps
$\mu_j:M\into \frak{g}^*$, $j=1,2,3$, we consider the same reduction
data as in Example~\ref{ex:hk},
$$
(\psi,\frak{h}=\frak{g}\oplus\frak{g}\oplus
\frak{g},\mu=(\mu_1,\mu_2,\mu_3)),
$$
but we will now use it to reduce the \gk\ structure $(\G,\J)$,
rather than the original hyper-K\"ahler structure.

Following Example \ref{ex:hk}, we know that $K^\G$ is given by
\eqref{eq:KG}, and each summand in that expression is invariant
under $I_j$, $j=1,2,3$. Note that we can write
$$
K^\G = (K^\G \cap C_+)\oplus(K^\G\cap C_-).
$$
If $X + g(X) \in K^{\G}\cap C_+$, then $I_1X -I_1^*g(X) = I_1X +
g(I_1X) \in K^\G$, and, similarly, $X-g(X) \in K^\G\cap C_-$ implies
that $I_2X+ I_2^*g(X)=I_2X-g(I_2X)\in K^\G$. Hence $\J K^{\G} =
K^{\G}$ and, by Theorem~\ref{theo:gk}, we may reduce the \gcs\
$(\G,\J)$ to $M^{red}=\mu^{-1}(0)/G$.

The bihermitian structure on $M_{red}$ associated with
$(\G^{red},\J^{red})$ can be described as follows: $\G$ restricted
to $K^{\G}$ is just the reduced metric $g^{red}$ obtained by
hyper-K\"ahler reduction as in Example \ref{ex:hk}, whereas the
restriction of $\J$ to $C_{\pm}\cap K^{\G}$ defines complex
structures via the projection to $TM^{red}$.  It is easily verified
that the restriction of $\J$ to $C_{+}\cap K^{\G}$ defines
$I_1^{red}$, and the restriction to $C_- \cap K^{\G}$ gives
$I_2^{red}$, where $(g^{red},I_1^{red},I_2^{red},I_3^{red})$ is the
reduced hyper-K\"ahler structure, as required.
\end{example}

\end{document}